\documentclass{llncs}   

\usepackage{amssymb}   
\usepackage[dvips]{graphicx}

\newcommand{\lex}{\hbox{\rm lex}}
\newcommand{\lpp}{\hbox{\rm lpp}}

\newcommand{\V}{\mathbb V}

\def\cocoa{\mbox{\rm C\kern-.13em o\kern-.07 em C\kern-.13em o\kern-.15em A}}

\begin{document}

\title{Automatic Discovery of Geometry Theorems Using Minimal Canonical Comprehensive Gr\"obner
Systems}
\titlerunning{Automatic Discovery with MCCGS}  
%
\author{Antonio Montes\footnote{Work partially supported by the Spanish Ministerio de
Ciencia y Tecnolog\'{\i}a under project MTM 2006-01267, and by the
Generalitat de Catalunya under project 2005 SGR 00692.} \inst{1}
\and Tom\'as Recio\footnote{Work partially supported by the
Spanish
  Ministerio de Educaci\'on y Ciencia under project GARACS,
  MTM2005-08690-C02-02.}
\inst{2} }
\authorrunning{A. Montes and T. Recio}   
%
\tocauthor{Antonio Montes (Universitat Polit\`ecnica de
Catalunya), Tom\'as Recio (Universidad de Cantabria)}
\institute{Dep. Matem\`atica Aplicada 2, Universitat Polit\`ecnica de Catalunya, Spain.\\
\email{antonio.montes@upc.edu}  \ \
\texttt{http://www-ma2.upc.edu/$\sim$montes} \\
\and Dep. Matem\'aticas, Estad{\'\i}stica y Computaci\'on, Universidad de Cantabria, Spain.\\
\email{tomas.recio@unican.es} \ \ \texttt{http://www.recio.tk}}

\maketitle              

\begin{abstract}
The main proposal in this paper is the merging of two techniques that
have been recently developed. On the one hand, we consider a
new approach for computing some specializable Gr\"obner basis, the so called Minimal Canonical
Comprehensive Gr\"obner Systems (MCCGS) that is -roughly speaking-  a
computational procedure yielding ``good" bases for ideals of polynomials
over a field, depending on several parameters, that specialize ``well", for instance,
regarding the number of solutions for the
given ideal, for different values of the parameters. The second
ingredient is related to automatic theorem discovery in elementary
geometry.  Automatic discovery aims to obtain complementary (equality and inequality type)
hypotheses for a (generally false) geometric statement to become
true. The paper shows how to use MCCGS for automatic discovering
of theorems and gives relevant examples.
\end{abstract}



\section{Introduction}
\subsection{Overview of Goals}
The main idea in this paper is that of merging two recent
techniques. On the one hand, we will consider a method (named {\it
MCCGS}, standing for {\it minimal canonical comprehensive
Gr\"obner systems}) \cite{MaMo06}, that is --roughly speaking--  a
computational approach yielding ``good" bases for  ideals of
polynomials over a field depending on several parameters, where
``good" means that the obtained bases should specialize (and
specialize ``well", for instance, regarding the number of
solutions for the given ideal) for different values of the
parameters.

Briefly, in order to understand what kind of problem {\it MCCGS}
addresses, let us consider the ideal  $(a x, x+y) K[a][x, y]$,
where $a$ is taken as a  parameter and $K$ is a field. Then it is clear that there
will be, for different values of $a=a_0 \in K$, essentially two different
types of bases for the specialized ideal $(a_0 x, x+y)
K[x, y]$.  In fact,  for $a_0=0$  we will get $(x+y)$ as a
Gr\"obner-basis (in short, a G-basis) for the specialized ideal; and for
any other rational value of $a$
such that $a=a_0 \neq 0$, we will get a G-basis with two elements,
$(x, y)$. Thus, the given G-basis $(a x, x+y) K[a, x, y]$ does not
specialize well to a G-basis of every specialized ideal. On the other
hand, let us consider $(a x-b) K[a,b][x]$, where $a, b$ are
taken as free parameters and $x$ is the only variable. Then, no
matter which rational values $a_{0}, b_{0}$  are assigned to $a,
b$, it happens that $\{a_{0} x-b_{0}\}$ remains a Gr\"obner basis
for the ideal $(a_{0} x-b_{0}) K[x]$. Still, there is a need for a
case-distinction  if we focus on the cardinal of the solutions for
the specialized ideal. Namely, for $a_0 \ne 0$ there is a unique
solution $x=-b_0/a_0$; for $a_0 = 0$ and $b_0 \ne 0$ there is no
solution at all;  and for $a_0 = b_0 = 0$ a solution can be any
value of $x$ (no restriction, one degree of freedom).

The goal of {\it MCCGS} is to describe, in a compact and canonical
form,  the discussion, depending on the different values of the
parameters specializing a given parametric system, of the different
basis for the resulting specialized systems and on their solutions.

The second ingredient of our contribution is about automatic theorem
discovery in elementary geometry.  Automatic discovery aims to
obtain complementary hypotheses for a (generally false) geometric
statement to become true. For instance, we can consider an
arbitrary triangle and the feet on each sides of the three
altitudes. These three feet give us another triangle, and now we
want to conclude that such triangle is equilateral. This is
generally false, but, under what extra hypotheses (of equality type) on the given
triangle will it become generally true?

\begin{figure}
\begin{center}
\includegraphics{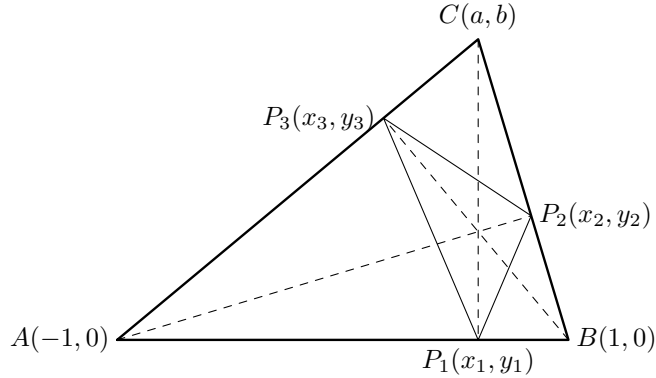}

\caption{\label{isotrianglefig} Orthic triangle}
\end{center}
\end{figure}

Finding, in an automatic way, the necessary and sufficient
conditions for this statement to become a theorem,  is the task of
automatic discovery. A protocol for automatic discovery is presented
in \cite{RV} and a detailed discussion of the method appears in \cite{DR}.
The protocol  proceeds requiring some computations (contraction, saturation, etc.)
about certain ideals built up from the given statement, but does not state any
preference about how to perform such computations (although the computed examples
in both papers rely on straightforward Gr\"obner bases computations for ideal elimination).

Our goal in this paper  is to show how we can improve the automatic discovery
of geometry theorems, by performing a {\it MCCGS}
procedure on an ideal built up from the given hypotheses
and theses, considering as parameters the free coordinates of some elements
of the geometric setting,

\subsection{Related work}
This idea has a close precedent in the work\footnote{But notice the authors of
\cite{CLLW}  already mention the paper of
Montes \cite{Mo02} as  a predecessor on this particular kind of
discussion of Gr\"obner basis with parameters.} of \cite{CLLW}, directly inspired by
\cite{K95} and, to a lesser extent, by \cite{Weis92}.
 In \cite{CLLW},  a parametric
radical membership test is presented for a mathematical construct
the authors introduce, called ``partitioned parametric Gr\"obner
basis" (PPGB).  Suppose we are given a statement $H:=\{h_1=0,
\dots, h_r=0\} \Rightarrow T:=\{g=0\}$, expressed in terms of
polynomial equations --usually over some computable field-- and their solutions over some extension field $K$
--that we can assume, in order to simplify the exposition, to be
algebraically closed. Roughly speaking, the method behind
\cite{CLLW} starts by computing the ``partitioned basis" (with
respect to a given subset of variables, here denoted by ${\bf u}$)
of an ideal $I\subseteq K[{\bf u}, {\bf x}, y]$, (for instance,
$I=(h_1({\bf u}, {\bf x})\dots h_r({\bf u}, {\bf x}), g({\bf u},
{\bf x})\, y-1)$),  ie.  a finite collection of couples
$(C_i,F_i)$, where the $C_i$'s are constructible sets described as
$\{c_1=0, \dots, c_m=0, q_1\neq 0, \dots , q_s \neq 0 \}$ on the
parameter space, and the $F_i$'s are some collections of
polynomials in $K[{\bf u}, {\bf x}, y]$. Moreover, it is required
(among other conditions) that the $C_i$'s conform a  partition of
the parameter space and, also, that for every element ${\bf u_0}$
in each $C_i$, the (reduced) G-basis of $(h_1({\bf u_0}, {\bf x})
\dots h_r,({\bf u_0}, {\bf x}),  g({\bf u_0}, {\bf x})\, y-1)$ is
precisely $F_i({\bf u_0}, {\bf x})$. It is well known (e.g.
\cite{K86} or \cite{Ch88}) that, in this context,  a statement
$\{h_1 = 0\dots h_r = 0\} \Rightarrow \{ g = 0\}$ is to be
considered true if $1 \in (h_1\dots h_r, g\,y-1)$; thus,  the
extra hypotheses that \cite{CLLW} proposes to add  for the
statement to become a theorem are precisely those expressed by any
of the $C_i$'s such that the corresponding $F_i =\{1\}$, since
this is the only case $F_i$ can specialize to $\{1\}$.

We must remark
that, simply testing for $1 \in (h_1\dots h_r, g\,y-1)$, as in the method above,  can
yield to theorems that hold just because the hypotheses are not
consistent (i.e. such that already $1 \in (h_1\dots h_r)$ ). This
cannot happen with our approach to automatic discovery: if a new statement is discovered, then the obtained
hypotheses will be necessarily consistent.

Although our approach stems from the same basic ideas, our contribution
differs from \cite{CLLW} in some respects: first, we
focus on automatic discovery, and not in automatic proving. Moreover,
we are able to specifically describe the capability and limitations of
the method (while in \cite{CLLW} it is only mentioned that, in the
reducible case,  their ``method \dots cannot determine if the conclusion
of the geometric statement is true on some components of the hypotheses").
Second, even for proving, the use of {\it MCCGS} provides not only the
specialization property (which is the key for the application of
partitioned parametric bases in  \cite{CLLW}) but also an automatic case
distinction, that allows a richer understanding of the underlying
geometry for the considered situation. In fact, it seems that the
partitioned parametric G-Basis (PPGB) algorithm from \cite{CLLW}
is close to the algorithm DISPGB considered in \cite{Mo02}, both
sharing that their output requires collecting by hand multiple
cases (and then having to manually express in
some simplified way the union of the corresponding conditions on
the parameters). Actually, the motivation for {\it MCCGS}  was,
precisely, improving DISPGB.

Our approach has also an evident connection (since \cite{Weis92} is the
common origin of all posterior developments on parametric Gr\"obner basis)
to the work of several members
of Prof. Weispfenning's group, regarding generic quantifier elimination
(Q. E.) and its application to
automatic theorem proving (as, for example, in \cite{DG}, \cite{DSW},
\cite{SS}, \cite{St}). In particular we remark the strong relation of
our work with that of \cite{DG}, that approaches theorem proving via
a restricted (generically valid) Q.E. method,  relying on generic
Gr\"obner systems computations. The set of restrictions $\Theta $
provided by this method, besides speeding up the Q.E. computations,
can be interpreted in the context of theorem proving, roughly speaking,
as a collection of new (sufficient) non-degeneracy conditions for an
statement to hold true.

Again, the difference between our contribution here and theirs is,
first, that we address problems requiring, in general, parameter
restrictions that go beyond  ``a conjunction $\Theta $ of negated
equations in the parameters" (\cite{DG}, first paragraph in
Section \ref{sec3}). That is, we deal with formulas that are
almost always false (see below for a more detailed explanation of
the difference between automatic derivation and automatic
discovery) and require non-negated (ie. equality) parameter
restrictions; they can not be directly approached via generic Q.E.
since our formulas are, quite often, generically false.  Moreover,
our approach is  limited to this specific kind of generically
false problems and we do not intend to provide a general method
for Q.E. A second difference is that, for our very particular kind
of problems, {\it MCCGS} formulates parameter restrictions in a
compact and canonical way, a goal that is not specifically
intended concerning the description of $ \Theta$ in \cite{DG}. For
these reasons we can not include performance comparisons to these
Q.E. methods and we do not consider relevant (although we provide
some basic information) giving hardware details, computing times,
etc. on the performance of our method running on the examples
described in the last section of this paper. We are not proposing
something better, but something different in a different context.

Next Section includes a short introduction to the basics on
automatic discovery, which could be of interest even for automatic
proving practitioners.  Section \ref{sec3} provides some
bibliographic references for the problem of the G-basis
specialization and summarizes  the main features of the  {\it
MCCGS} algorithm, including an example of its output. Section
\ref{sec4} describes the application of {\it MCCGS} to automatic
discovery, while Section 5 works in detail a collection of curious
examples, including the solution of a pastime from {\it Le Monde}
and the simpler solution (via this new method) of one example also
solved by a more traditional method.

\section{A digest on automatic discovery}

Although less popular than automatic proving, automatic discovery of
elementary geometry theorems is not new.  It can be traced back to
the work of Chou (see \cite{Ch84}, \cite{Ch87} and \cite{ChG90}),
regarding the ``automatic derivation of formulas", a particular
variant of automatic discovery where  the goal  consists in
deriving results that always occur under some given hypotheses
but that can be formulated in terms of some specific
set of variables (such as expressing the area of a triangle in
terms of the lengths of its sides).  Finding the geometric locus
of a point  defined through some geometric constraints (say,
finding the locus of a point when its projection on the three
sides of a given triangle form a triangle of given constant area \cite{Ch88},
Example 5.8) can be considered as another  variant of this
``automatic derivation" approach.

Although ``automatic derivation" (or locus finding) aims to discover
some new geometric statements (without modifying the given hypotheses),
it is not exactly the same as ``automatic discovery" (in the sense we
have presented it in the previous section), that searches for
complementary hypotheses for a (generally false) geometric
statement to become true (such as stating that the three feet of
the altitudes for a given triangle form an equilateral
triangle and finding what kind of triangles verify it).
Again, automatic discovery in this precise sense appears in the
early work of Chou (whose thesis \cite{Ch85} deals with
``Proving and discovering theorems in elementary geometries using
Wu's method")  and Kapur \cite{K89} (where
it is explicitly stated that ``\dots the objective here is to find
the missing hypotheses so that a given conclusion follows from a
given incomplete set of hypotheses\dots").

Further specific contributions to automatic discovery appear in
\cite{Wa98}, \cite{R98} (a book written in Spanish for
secondary education teachers, with circa one hundred pages devoted
to this topic and with many worked out examples), \cite{RV},
\cite{Koepf} or \cite{CW}.  Examples of automatic derivation,
locus finding and discovery, achieved through a specific
software named {\it GDI} (the initials of {\it
Geometr{\'\i}a Din\'amica Inteligente}), of Botana-Valcarce,
appear in \cite{BR2005} or \cite{RB} (and the references thereof),
such as the automatic
derivation of the thesis for the celebrated Maclane $8_3$-Theorem,
or the automatic answer to some items on a test posed by
Richard \cite{Richard}, on proof strategies in mathematics
courses, for students 14-16 years old.

The simple idea behind the different approaches
is\footnote{Already present in the well known book of \cite{Ch88},
page 72:  ``\dots The method developed here can be modified for the
purpose of finding new geometry theorems\dots Suppose that we are
trying to prove a theorem\dots and the final remainder\dots $R_0$
is nonzero. If we add a new hypotheses $R_0=0$, then we have a
theorem\dots ".  Here Chou proposes adding as new hypotheses the
pseudoremainder of the thesis by the ideal of hypotheses, a
mathematical object which should be zero if the theorem was
generally true.}, essentially,
that of adding the conjectural theses to the collection of
hypotheses, and then deriving, from this new ideal of theses plus
hypotheses, some new constraints in terms of the free parameters
ruling the geometric situation. For a toy example, consider that
$x-a=0$ is the only hypothesis, that the set of points $(x, a)$
in this hypothesis variety is determined by the value of the
parameter $a$, and that
$x=0$ is the (generally false) thesis. Then we add the thesis to
the hypothesis, getting the new ideal $(x-a, x)$, and we observe
that the elimination of $x$ in this ideal yields the constraint
$a=0$, which is
indeed the extra hypothesis we have to add to the given one $x-a=0$,
in order to have a correct statement $[x-a=0 \wedge a=0] \Rightarrow [x=0]$.

 With this simple idea as starting point\footnote{Indeed, things are
 not so trivial. Consider, for instance,  $H \Rightarrow T$, where
 $H= (a+1)(a+2)(b+1) \subset K[ a, b, c ]$ and
 $T=(a+b+1, c)\subset K[ a, b, c ]$. Take as parameters $U=\{b, c\}$,
 a set of ${\rm dim}(H)$-variables, independent over $H$.
 Then the elimination of the remaining variables over $H+T$ yields
 $H'= (c, b^3-b)$.
 But $H+H'=(a+1, b, c) \cap (a+2, b, c)\cap (a+1, b-1, c) \cap (a+2, b-1, c) \cap (b+1, c)$
 does not imply $T$, even if we add some non-degeneracy conditions
 expressed in terms of the free parameters $U$, since $T$ vanishes over
 some components,  such as  $(a+2, b-1, c)$ (and does not vanish over
 some other ones, such as $(a+1, b-1, c)$). },  an elaborated discovery
 procedure, with several non trivial examples, is presented in \cite{RV}.
It has been recently revised in \cite{BDR} and \cite{DR},  showing
that, in some precise sense, the idea of considering $H + T$ for
discovering is intrinsically unique (see Section \ref{sec4} for a
short introduction, leading to the use of {\it MCCGS} in this
context).

\section{Overwiew on the {\it MCCGS} algorithm}\label{sec3}

As mentioned in the introduction, specializing the basis of an
ideal with parameters does not yield, in general, a basis of the
specialized ideal.

This phenomenon --in the context of Gr\"obner basis-- has been
known for over fifteen years now, yielding to a rich variety of
attempts towards a solution (we refer the interested reader to the
bibliographic references in \cite{MaMo06} or in \cite{Wib06}).
Finding a specializable basis (ie. providing a single basis that
collects all possible bases, together with the corresponding
relations among the parameters) is --more or less--  the task of
the different comprehensive G-Basis proposals. Although  the first
global solution was that of Weispfenning, as early as 1992 (see
\cite{Weis92}), the topic is quite active nowadays, as exemplified
in the above quoted recent papers.  The {\it MCCGS} procedure,
that is, computing the {\it minimal canonical comprehensive
Gr\"obner system} of a given parametric ideal,  is one of the
approaches we are interested in.  Let us describe briefly the
goals and output of the {\it MCCGS} algorithm.

Given a parametric polynomial system of equations over some
computable field, such as the rational numbers, our interest
focuses on discussing the type of solutions over some
algebraically closed extension, such as the complex numbers,
depending on the values of the parameters. Let ${\bf
x}=(x_1,\dots,x_n)$ be the set of variables, ${\bf
u}=(u_1,\dots,u_m)$ the set of parameters and $I\subset K[{\bf
u}][{\bf x}]$ the parametric ideal we want to discuss, where, in
order to simplify the exposition, a single field $K$,
algebraically closed,  is considered both for the coefficients and
the solutions. We want to study how the  solutions over $K^n$ of
the equation system defined by $I$ vary when we specialize the
values of the parameters ${\bf u}$ to concrete values ${\bf u_0}
\in K$. Denote by $A=K[{\bf u}]$, and by $\sigma_{{\bf u_0}}:
A[{\bf x}] \rightarrow K[{\bf x}]$ the homomorphism corresponding
to the specialization (substitution of ${\bf u}$ by some ${\bf
{u_0}} \in K$).

A Gr\"obner System $GS(I,\succ_{{\bf x}})$ of the ideal $I\subset
A[{\bf x}]$  wrt (with respect to) the termorder $\succ_{{\bf x}}$
is a set of pairs $(S_i, B_i)$, where each couple consists of a
constructible set (called segment) and of a collection of
polynomials, such that
\[
\begin{array}{lcl}
\hbox{\rm GS}(I,\succ_{{\bf x}})&=& \{ (S_i,B_i)\ :\ 1 \le i \le
s,\ S_i \subset K^m,\ B_i
\subset A[{\bf x}],\ \bigcup_i S_i=K^m,  \\
&&  \forall {\bf u_0} \in S_i,\  \sigma_{{\bf u_0}}(B_i)  \hbox{ is a
Gr\"obner basis of } \sigma_{{\bf u_0}}(I) \hbox{ wrt } \succ_{{\bf x}} \}.
\end{array}
\]

The algorithm {\it MCCGS} (Minimal Canonical Comprehensive
Gr\"obner System) \cite{Mo06},\cite{MaMo06} of the ideal $I\subset
A[\bf{x}]$ wrt the monomial order $\succ_{\bf{x}}$ for the
variables, builds up the unique Gr\"{o}bner System having the
following properties:

\begin{enumerate}
  \item  The segments $S_i$ form a partition ${\mathcal S}=\{S_1,\dots,S_s\}$ of the parameter space
  $K^m$.

 \item The polynomials in $B_i$ are normalized to
have content 1 wrt $\bf{x}$ over $K[\bf{u}]$ (in order to work
with polynomials instead of with rational functions). The $B_i$
specialize to the reduced Gr\"obner basis of
$\sigma_{{\bf{u_0}}}(I)$, keeping the same  $lpp$ (leading power
products set) for each ${\bf u_0} \in S_i$, i.e. the leading
coefficients are different from zero on every point of $S_i$.
\footnote{The polynomials in the $B_i$'s are not faithful (they do
not belong to $I$), as they are reduced wrt to the null conditions
in $S_i$. By abuse of language we call them {\it reduced bases}
 (i.e. not-faithful, in the terminology of Weispfenning).}.  Thus a
concrete set of $lpp$ can be associated to a given $S_i$. Often it
exists a unique segment corresponding to each particular $lpp$,
although in some cases several such segments can occur. In any
case, when a segment with the reduced basis [1] exists, then it is
unique. When two segments $S_i$, $S_j$ share the same $lpp$, then
there is not a common reduced basis  $B$ specializing to both
$B_i, B_j$\footnote{M. Wibmer~\cite{Wib06} has proved that for
homogeneous ideals in the projective space there is at most a
unique reduced basis and segment corresponding to a given $lpp$.}.

Moreover, there exists a unique  segment $S_1$ (called the {\it generic segment}), containing  a
Zariski-open set,  whose associated basis $B_1$  is called the {\it generic basis} and coincides with the
Gr\"obner basis of $I$ considered in $K({\bf u})[{\bf x}]$
conveniently normalized without denominators and content 1 wrt
${\bf x}$.

 \item  The partition
 ${\mathcal S}$ is canonical (unique for a given $I$ and monomial
 order).
  \item The partition is  minimal, in the sense it
  does not exists another partition having property 2
  with less sets $S_i$.
  \item The segments $S_i$ are described in a canonical form.
\end{enumerate}

As it is known, the $lpp$ of the reduced Gr\"obner basis of an
ideal determine the cardinal or dimension of the solution set over
an algebraically closed field. This makes the {\it MCCGS}
algorithm very useful for applications as it identifies
canonically the different kind of solutions  for every value of
the parameters. This is particularly suitable for automatic
theorem proving and automatic theorem predicting,  as we will show
in the following sections.

Let us give an example of the output of {\it MCCGS}.

\begin{example} \label{singpointsconic}
Consider the system described by the following parametric ideal
(here the parameters are $a, b, c, d$):
\[I=( x^2+b y^2+2 c x y+2 d x,2 x+2 c y+2 d,2 b y+2 c x
), \] arising in the context of finding all possible singular
conics and their singularities.  Calling to the Maple
implementation of {\it MCCGS} yields a graphical and an algebraic
output. The graphical output is shown in Figure \ref{conictree}.
It contains the basic information that is to be read as follows.
At the root there is the given ideal (in red). The second level
(also in red) contains the $lpp$ of the bases of the three
different possible cases. These are $[1]$, corresponding to the no
solution (no singular points) case; $[x,y]$, corresponding to the
one solution (one singular point) case;  and $[x]$, corresponding
to the case of one dimensional solution (ie. when the conic is a
double line). Below each case there is a subtree (in blue)
describing the corresponding $S_i$, with the following
conventions:

\begin{itemize}
\item at the nodes there are  ideals of $K[{\bf u}]$, prime in the field
of definition (generated over the prime field by the coefficients of a
reduced G-Basis) of the given ideal $I \subset A[{\bf x}]$

\item a descending edge means the set theoretic ``difference" of the set
defined by the node above minus the set defined at the node below,

\item nodes at the same level, hanging from a common node, are to be
interpreted as yielding the set theoretic ``union" of the
corresponding sets; they form the irredundant prime decomposition
of a radical ideal of $K[{\bf u}]$.
\item every branch contains a strictly ascending chain of prime ideals.
\end{itemize}

 So, in the example above, the three  cases, their $lpp$ and the corresponding
$S_i$'s are to be read as shown in the following table:

\begin{center}
\begin{tabular}{||c|l|l||}
 \hline \hline
 $\lpp$ & Basis $B_i$ & Description of $S_i$ \\
 \hline
 $[1]$ & $[1]$ & $K^3 \setminus \left( \left( \V(b)
\setminus (\V(c,b) \setminus \V(d,c,b)) \right) \cup \V(d) \right)
 $ \\
 \hline
 $[y,x]$ & $[2 c y+d, x]$ & $\left(\V(b) \setminus \V(c,b) \right) \cup  \left( \V(d) \setminus
 \V(d,b-c^2)\right)
  $ \\
  \hline
  $[x]$ &  $[x+cy]$ & $ \V(d,b-c^2) $ \\
 \hline\hline
\end{tabular}
\end{center}\vspace{3mm}

We remark that the $B_i$'s do not appear in the Figure
\ref{conictree}, since --in order to simplify the display-- the
complete bases are only given by the algebraic output of {\it
MCCGS} and are not shown by the graphic output.

\begin{figure}
\begin{center}
\includegraphics{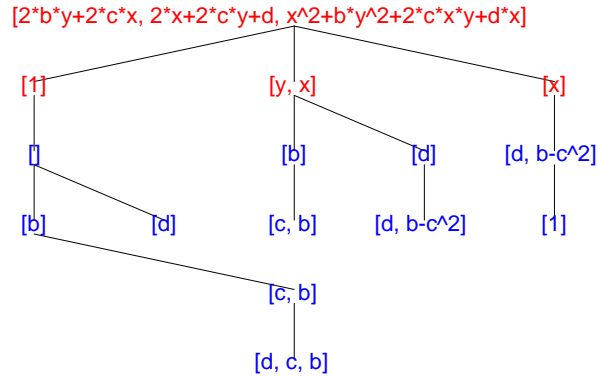}
\caption{\label{conictree} {\it MCCGS} for the singular points of
a conic}
\end{center}
\end{figure}
\end{example}

\section{Using {\it MCCGS} for automatic theorem
discovering}\label{sec4}

Once we have briefly described the context for {\it MCCGS} and for
automatic discovery, we are prepared to describe the basic idea in
this paper. We can say that our goal  is to show how performing a
{\it MCCGS} procedure can improve the automatic discovery of
geometry theorems.

Example \ref{singpointsconic} can be seen as a very simple example
of theorem discovering. We could formulate the statement {\it a
conic has one singular point} and try to find the conditions for
the statement to be true. Without loss of generality we express
the equation of the conic and its partial derivatives as
\[I=( x^2+b y^2+2 c x y+2 d x, 2 x+2 c y+2 d, 2 b y+2 c x
), \] and search for the values of the parameters where this
system has a single solution. As shown above, we have found that
the statement is true if and only if $\{b=0, c \ne 0\}$ or if
$\{d=0, b-c^2 \ne 0\}$, since in the first segment of the table there is no
solution ($B_1=(1)$), while the third segment yields a 1-dimensional
set of solutions.

In general, let $H \Rightarrow T$ be a statement expressed in
terms of polynomial equations,  where the ideals $H, \; T
\subseteq K[x_1, \dots, x_n]$  will be the corresponding
hypotheses ideal and theses ideal (both, possibly, with several
generators).  In this context \cite{DR} sets, as discovery goal,
finding a couple of subsets of variables $U \supseteq U'$, with
$X\supseteq U \supseteq U'$, and a couple of ideals $R'\subset
K[U], R''\subset K[U']$, so that the following properties hold for
the associated algebraic varieties (over $K^n$, with $K$
algebraically closed):

\begin{enumerate}
\item $\V(H+R'^e ) \backslash \V(R''^e) \subseteq \V(T )$ (where the $^e$ stands for
the extension of the ideal from its defining ring, say, $K[U]$ or $K[U']$,  to $K[X]$);
\item $\V(H+T ) \subseteq \V(R'^e)$;
\item $\V(H+R'^e) \backslash \V(R''^e) \neq \emptyset$.
\end{enumerate}

The rationale behind such a definition is that such a couple $(R', R'')$ is supposed to provide
\begin{itemize}
\item some necessary (as expressed by item $2)$ above)
\item and sufficient  (as expressed by item $1)$  above)
\end{itemize}
\noindent non trivial (as expressed by item $3)$ ) complementary conditions of
equality kind (given by  $R'$ ) and of  non-degeneracy type (given by the negation
of $R''$) for the given theses to hold under the
given hypotheses.

Then it is shown in \cite{DR} that, for a given couple of subsets
of variables $U \supseteq U'$, with  $X\supseteq U \supseteq U'$,
there is a couple of ideals $R'\subset K[U], R''\subset K[U']$
verifying properties $1),  2)$ and  $3)$ above if and only if  the
couple of ideals $H' = (H +T) \cap K[U]$ and\footnote{Let  $I, J$
be ideals of $K[X]$. Recall that $I:J=\{x, xJ\subset I\}$. Then,
the {\em saturation of $I$ by $J$} is defined as
$I:J^\infty=\cup_{n} (I:J^n)$, cf. \cite{KR} .} $H'' =
((H+H'^e):T^{\infty})\cap K[U']$ also verify these three
conditions. Moreover, Theorem 2 in \cite{DR} shows these
conditions hold  if and only if $1\not\in (H')^{c}:H''^{\infty}$
(equivalently, iff $H'' \not\subseteq \sqrt{(H')^{c}}$), where
$^c$ stands for the contraction ideal, so there is an algorithmic
way of solving the posed discovery problem for a given statement
and choice of variables.

  Now we remark the following:
\begin{proposition}\label{discovparam}
If there is  a couple $R', R''$ verifying the above conditions, then
$$\V(H+R'^e) \backslash \V(R''^e) =  \V(H+H'^e) \backslash \V(R''^e)$$
\end{proposition}
\begin{proof}
First notice that $\V(H+H'^e)\subseteq \V(H+R'^e)$, since,  by
property $2)$, $\V(H+T) \subseteq \V(R'^e)$, thus $R'^e \subseteq
\sqrt{H+T }$ and so $R'=R'^{ec} \subseteq \sqrt{H+T }^c =
\sqrt{H'}$ and this implies that $\V(H'^e) \subseteq \V(R'^e)$.

Moreover, we have also that $\V(H+R'^e) \backslash \V(R''^e)
\subseteq \V(T) \cap \V(H)$, by property $1)$, and $\V(T) \cap
\V(H) \subseteq \V(H'^e)$, where the last inclusion follows from
the definition of $H'^e$.  We conclude that $\V(H+R'^e) \backslash
\V(R''^e) \subseteq \V(H+H'^e) \backslash \V(R''^e)$.
\end{proof}

This means that the search for candidates $R'$ for complementary
hypotheses of equality type, can be reduced to computing $\V(H')$.
This is, precisely, the (Zariski closure of the) projection, over
the parameter space of the $U$-variables, of $\V(H + T)$,  and
this can be computed through {\it MCCGS}, providing as well some
other useful information (as in Corollary \ref{cor4}).

\begin{proposition}\label{discovcomp}
The projection of $\V(H+T)$ over the $U$-variables can be computed
by performing a {\it MCCGS} for $I = H+T$ and $X \supset U$,
discarding, if it exists, the unique segment $S_i$'s with $B_i$
equal to $1$ and keeping the remaining $S_i$'s.
\end{proposition}

\begin{proof}
Since the segments of a {\it MCCGS}  partition  the parameter
space $U$, it is enough to show that a point $(u_0)$ is not in the
projection if and only if it belongs to the $S_i$ with associated
$B_i = 1$. Now we recall that a reduced Gr\" obner basis is 1 if
and only if the corresponding ideal is $(1)$. Then, the ideal $H +
T$ specialized at $u_0$ will be $(1)$ if and only if its reduced
G-basis is 1. Since we work over an algebraically closed field,
this is the only case the system $H + T$, specialized at $u_0$,
has no solution, ie. $u_0$ is not in the projection of $\V(H +
T)$. But, by construction, a  $B_i$ specializes to 1 if and only
if $B_i = 1$ (since the specialization must be a reduced G-basis
and has the same $lpp$ as $B_i$).
\end{proof}

\begin{corollary}\label{cor4} The union of these $S_i$'s with associated  $B_i \neq 1$
(ie. the complement of the only possible segment with $B_i =1$)
partitions the projection of $\V(H + T)$;  that is, it holds $H
\wedge T \Rightarrow \{ \cup S_i \}$. Thus,  the union of these
$S_i$'s provide complementary necessary conditions for the theses
$T$ to hold over $H$.
\end{corollary}

We will see below (Remark \ref{rem9}) that, when the given
statement does not hold over any geometrically meaningful
component of the hypotheses variety -- ie. in the automatic
discovery situation-- the segment with $B_i =1$ is the generic
one, so its complement provides necessary conditions for the
theses $T$ to hold over $H$.

 Next we must study if some of these $S_i$'s provide sufficient conditions, analyzing
 the behavior of each statement $H \wedge S_i \Rightarrow T$,  for every  segment
 $S_i$ with $\lpp \neq 1$. Some --perhaps all,  perhaps none-- of them could be true.
 Remark that, anyway, $H \wedge S_i \neq \emptyset$, since the associated basis in not $1$.
 Remark, also,  that {\it MCCGS} allows to obtain supplementary conditions $S_i$ of
 the more general form (not every constructible set is the difference of two closed
 sets of the form $\V(H + R'^e) \backslash \V(R''^e)$, as in the previous approach).

 There are some special easy cases, as shown in the next result.

\begin{corollary}\label{segvars}
For every segment $S_i$ such that the corresponding $lpp$ of the
associated basis is, precisely, the collection of variables
$\{x_1, \dots, x_n\}$, we have that $\V(H) \cap S_i \subseteq
\V(T)$, ie. $H \wedge S_i \Rightarrow T$ holds, and $S_i$ provides
sufficient conditions for $T$ to hold over $H$.
\end{corollary}
\begin{proof}
In fact, the condition on the associated $lpp$ means that for
every $u_0$ in $S_i$, the system $H(u_0, x)=0, T(u_0, x)=0$ has a
unique solution, and it belongs to $\V(T)$. Thus $\V(H) \cap S_i
\subseteq \V(T)$.
\end{proof}

Otherwise, we should analyze, for each $i$ with $S_i$ involved in
the projection of $\V(H+T)$, the validity of $H \wedge S_i
\Rightarrow T$.  This is a straightforward  ``automatic proving"
step, and not of ``automatic discovery",  since adding again $T$
to the collection of hypotheses $H \wedge S_i $ will not change
the situation, as the projection of $\V(H) \cap \V(T) \cap  S_i$
equals the projection of $\V(H) \cap  S_i$, both being $S_i$.

Yet, {\it MCCGS} can provide a method for checking the truth of such statement
$H \wedge S_i \Rightarrow T$.  As it is well known, we can reformulate the hypotheses
$H \wedge S_i$ as a collection of equality hypotheses $H$, since $S_i$ is constructible
and, then, the union of intersections of closed and open sets (in the Zariski topology).
And open sets can be expressed through equalities by means of saturation techniques
(such as $x \neq 0 \Leftrightarrow x y -1 =0$, etc.). So let us state the following
propositions (adapting to the {\it MCCGS} context some results from \cite{RV},
\cite{DR}) in all generality.

\begin{proposition}\label{prop6}
Let $H \Rightarrow T$ be a statement and let $U$ be a collection
of variables independent for  $H$. Then $T$ vanishes identically
on all the components of $H$ where $U$ remain independent if and
only if, performing a {\it MCCGS} for $\{H, T z -1\}$ with respect
to $U$, the generic basis is 1.
\end{proposition}

\begin{proof}
Notice the stated condition on the segments of the {\it MCCGS} is
equivalent to the fact that the contraction $(H, T z -1) \cap K[U]
\neq (0)$. In fact, this contraction is zero if and only if the
projection of $\V(H + (T z -1))$ contains an open set. And this is
equivalent to the fact that the generic segment has $\lpp \neq 1$.

Now, if $(H, T z -1) \cap K[U] \neq (0)$, take some $0 \neq g \in
(H, T z -1) \cap K[U]$. Remark that, by construction,  $g\,T = 0$
over $\V(H)$.  If $T \neq 0$ at some point over some component of
$\V(H)$, then $g = 0$ over such component; so it cannot  be a
component where the  $U$ are independent, since $g \in K[U]$.

Conversely, if $T$ vanishes identically over all the independent
components, then we can compute an element $g \in K[U]$ vanishing
over the remaining components (because $U'$ is dependent over
them). So $g\,T$ vanishes all over $\V(H)$, and thus $(H, T z -1)
\cap K[U'] \neq (0)$.
\end{proof}

\begin{remark}\label{segs1}
In fact, as in \cite{CLLW},  it is easy to show that the segment
with associated $lpp$ equal to 1 provides complementary sufficient
conditions for $H \Rightarrow T$ to hold. In fact, for every ${\bf
u_0}$ in such segment, $\V(H({\bf u_0}, {\bf x}), T({\bf u_0},
{\bf x}) z-1) = \emptyset$, so $\V(H({\bf u_0}, {\bf x}) \subseteq
T({\bf u_0}, {\bf x})$. But it can happen there is no such
segment.
\end{remark}

\begin{proposition}\label{genseg1}
Let $H \Rightarrow T$ be a statement and let $U$ be a collection
of variables independent for  $H$ and of dimension equal to ${\rm
dim}(H)$. Then $T$ vanishes identically on some components of $H$
where $U$ remains independent if and only if, performing a {\it
MCCGS} for $\{H, T\}$ with respect to $U$, the reduced basis of
the generic segment is different from 1.
\end{proposition}

\begin{proof}
As above, the stated condition on the segments of the {\it MCCGS} is equivalent
to the fact that the contraction $(H, T ) \cap K[U']  = (0)$.

Now, if  $T$ does not vanish identically over any component of
$\V(H)$ independent over $U'$, the projection of $\V(H, T)$ over
$U$ will be a proper closed subset (since the dimension of the
projection is less or equal than the dimension of the components
of $\V(H)$ contained in $\V(T)$, the  maximum dimension of all
components of $\V(H)$ equals the maximum dimension of the
independent components, and the dimension of the $U$-space equals
the maximum dimension of the components of $\V(H)$). This
contradicts the assumption $(H, T ) \cap K[U] =  (0)$, which
implies the closure of the projection is the whole $U$-space.

Conversely, if $T$ vanishes identically over some independent
component (say, $C$)  and $(H, T ) \cap K[U] \neq  (0)$, then we
can choose an element $0 \neq g \in (H, T ) \cap K[U]$. This
element vanishes over any component of $\V(H)$ where $T$ vanishes,
in particular over $C$, contradicting its independence over $U$.
\end{proof}
\begin{remark}\label{rem9}
The last proposition  can be also read in a different way: $T$
does not vanish identically on any independent component of $H$ if
and only if the reduced basis of the generic segment is  1.
\end{remark}

\begin{corollary}\label{cor10}
Let $H \Rightarrow T$ be a statement and let $U$ be a collection
of variables independent for  $H$ and of dimension equal to ${\rm
dim}(H)$. Then $T$ vanishes identically on some components of $H$
where $U$ remains independent and also  $T$ does not vanish
identically on some other components of $H$ where $U$ remains
independent if and only if
\begin{itemize}
\item performing a {\it MCCGS} for $\{H, T z -1\}$ with respect to
$U$, the generic segment dos not have reduced basis 1,  and

\item performing a {\it MCCGS} for $\{H, T\}$ with respect to $U$, the reduced basis of
the generic segment is also different from 1.
\end{itemize}
\end{corollary}

In conclusion, using {\it MCCGS} one can determine, for a given
statement, whether it is  generally true (over all independent
components, using Proposition \ref{prop6}, generally false (over
all independent components, using Remark \ref{rem9}), or partially
true and false (using Corollary \ref{cor10}). Let us call this
last situation the ``undecidable" case.

In fact, unfortunately, in this circumstance it is not possible,
using only data on the $U$ variables, to determine the components
of $H$ where $T$ vanishes identically. Consider $H= b(b+1)
\subseteq K[a,b]$, $T= (b)$ and take $U= \{a\}$. Here the
projection of $\V(H, T z-1)$ over the $U$-variables is the whole
$a$-line, so does not have any segment with $lpp$ equal to 1, and
we know the thesis does not hold over all independent components.
Moreover the projection of $\V(H + T)$ over the $U$-space is again
the whole $a$-line, so there is no segment with $lpp$ 1, and we
can conclude $T$ holds over a component, but there is no way of
separating the component $b=0$, by manipulating $H, T$ in terms of
polynomials in the variable $a$.

This discussion applies to the situation described above, when
considering statements $H \wedge S_i \Rightarrow T$,  where
segment $S_i$ belongs to a {\it MCCGS} for $\{H, T\}$ with respect
to a collection $U$ of variables and has $\lpp \neq 1$.  Let $HH$
be the reformulation of $H \wedge S_i$ in terms of equalities and
let (if possible) $U' \subseteq U$ be a new collection of
variables,  such that they are independent for  $HH$ and of
dimension equal to ${\rm dim}(HH)$.

Then, as remarked above,  $HH \Rightarrow T$ will be true on the segment
$SS_i$ of a {\it MCCGS} with respect to $HH, T z-1$, with $lpp$ 1. If it is
is an open segment, then the statement  $H \wedge S_i \Rightarrow T$
will be generally true (over all the components independent over $U'$).
If it is not open segment, but there is at least one such segment,
the statement will hold true under the new restrictions.

But if there is no segment at all with $\lpp = 1$, then and only
then we are in the undecidable case. In fact,  over all points in
the $U'$-projection of $\V(HH, T z -1)$ we will have points  of
$\V(H)$ not in $\V(T)$ (because all the segments will have  $\lpp
\neq 1$ in the {\it MCCGS} for $HH, T z -1$) and also points of
$\V(H)$ and $\V(T)$ (since we are also in the projection of $S_i$
over $U'$, and $S_i$ corresponds to a segment of $\lpp \neq 1$ for
a {\it MCCGS} with respect to $H, T$).

In this case, since the projection over $U'$ of $\V(HH, T)$  will
be same as the projection of $\V(HH)$ (both being equal to the
projection of $S_i$), it is of no use to go further with a new
discovery procedure, computing a {\it MCCSG} for $HH, T$ over
$U'$. We know beforehand that all its segments will have $\lpp
\neq 1$, since over any point in the projection of $S_i$ there
will be always points on $\V(HH) \cap \V(T)$, confirming, again,
that we are in the undecidable situation.

\section{Examples}

Let us see how this works in a collection of examples, where we
have just detailed the discovery step (ie. computing the {\it
MCCGS} of $\{H, T\}$ with respect to a collection of maximal
independent variables for $H$,  and then collecting the
potentially true statements $H \wedge S_i \Rightarrow T$, where
segment $S_i$ has $\lpp \neq 1$)  in the procedure outlined in the
previous Section.  That is, we have not included here the formal
automatic verification in each case that the newly found
hypotheses actually lead to a true statement (the ``proving
step").

\begin{example}\label{example2} (See also \cite{DR}).
Next, we will develop the above introduced notions considering a
statement from   \cite{Ch88} (Example 91 in his book), suitably
adapted to the discovery framework. The example here is taken from
\cite{DR}.

   Let us consider as given data a
circle and two diametral opposed points on it (say, take a circle
centered at $(1,0)$ with radius $1$, and let $C=(0,0), D=(2, 0)$
the two ends of a diameter), plus an arbitrary point $A=(u_1,
u_2)$. See Figure \ref{ex4fig}. Then trace a tangent from $A$ to
the circle and let $E=(x_1,x_2)$ be the tangency point. Let
$F=(x_3,x_4)$ be the intersection of $DE$ and $CA$. Then we claim
that $AE=AF$. Moreover, in order to be able to define the lines
$DE$, $CA$, we require, as hypotheses,  that $D\neq E$ (ie. $u_1
\neq 2$) and that $C\neq A$ (ie. $u_1 \neq 0$ or $u_2 \neq 0$).

\begin{figure}
\begin{center}
\includegraphics{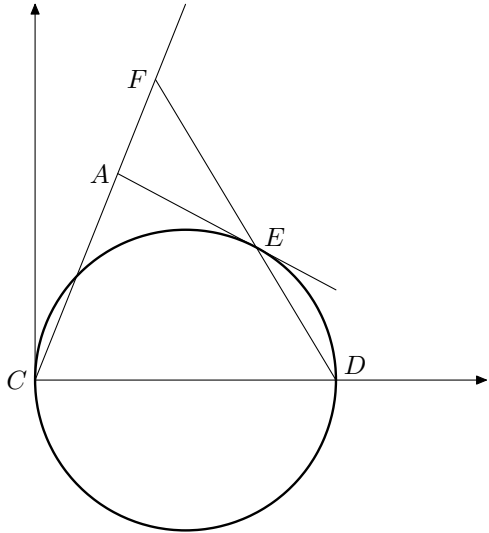}

\caption{\label{ex4fig} Problem of Example \ref{example2}}
\end{center}
\end{figure}

 Now, using \cocoa \ \cite{cocoa} and its package TP (for Theorem
Proving),
we translate
the given situation as follows
\begin{verbatim}
Alias  TP   := $contrib/thmproving;

Use R::=Q[x[1..4],u[1..2]];

A:=[u[1],u[2]];
E:=[x[1],x[2]];
D:=[2,0];
F:=[x[3],x[4]];
C:=[0,0];

Ip1:=TP.Perpendicular([E,A],[E,[1,0]]);
Ip2:=TP.LenSquare([E,[1,0]])-1;
Ip3:=TP.Collinear([0,0],A,F);
Ip4:=TP.Collinear(D,E,F);

H:=Saturation(Ideal(Ip1,Ip2,Ip3,Ip4),Ideal(u[1]-2)*
                                  Ideal(u[1], u[2]));

T:=Ideal(TP.LenSquare([A,E])-TP.LenSquare([A,F]));
\end{verbatim}
\noindent where $T$ is the thesis and $H$ describes the hypothesis
ideal. Notice that $Ip1$ expresses that the segments $[E, A], [E,
(1,0)]$ are perpendicular; $Ip2$ states that  the square of the
length of $[E, (1, 0)]$ is 1 ( so $Ip1, Ip2$ imply  $E$ is the
tangency point from $A$ );  and the next two hypotheses express
that the corresponding three points are collinear. The hypothesis
ideal $H$ is here constructed by using the $saturation$ command,
since it is a standard way of stating  that the hypothesis variety
is the (Zariski) closure of the set defined by all the conditions
$Ip[i]=0,  i=1\dots 4$ minus the union $\{u[1]=2\} \cup \{u[1]=0,
u[2]=0\}$, as declared in the formulation of this example (but we
refer  to \cite{DR} for a discussion on the two possible ways of
introducing inequalities as hypotheses). Finally, the thesis
expresses that the two segments $[AE], [AF]$ have equal non
oriented length.

Now, let us use in this, clearly false,  statement  the approach
of \cite{RV} or \cite {DR} to discovery. First we check that the
statement $H \rightarrow T$ is not algebraically true in any
conceivable way. For instance, it turns that
\begin{verbatim}
Saturation(H, Saturation(H,T));
Ideal(1)
-------------------------------
\end{verbatim}

\noindent and this computation shows that all possible non-degeneracy
conditions (those polynomials $p(\bf{u}, \bf {x})$ that could be added
to the hypotheses as conditions of the kind $p(\bf{u},\bf{x}) \neq 0)$
lie in the hypotheses ideal, yielding, therefore to an empty set of
conditions of the kind $p \neq 0 \wedge p =0$.  This implies, in
particular, that the same negative result would be obtained if we
restrict the computations to some subset of variables, since the thesis
does not vanish on any irreducible component of the  hypotheses variety.

Thus we must switch on to the discovery protocol, checking before hand
that  $u[1], u[2]$ actually is a (maximal) set of independent variables --the parameters--
for our construction:

\begin{verbatim}
Dim(R/H);
2
-------------------------------

Elim([x[1],x[2],x[3],x[4]],H);
Ideal(0)
-------------------------------
\end{verbatim}

\noindent Then we add the thesis to the hypotheses ideal and we
eliminate all variables except $u[1], u[2]$

\begin{verbatim}
H':=Elim([x[1],x[2],x[3],x[4]],H+T);
H';
Ideal(-1/2u[1]^5 - 1/2u[1]^3u[2]^2 + u[1]^4)
-------------------------------
Factor(-1/2u[1]^5 - 1/2u[1]^3u[2]^2 + u[1]^4);
[[u[1]^2 + u[2]^2 - 2u[1], 1], [u[1], 3], [-1/2, 1]]
-------------------------------
\end{verbatim}
\noindent yielding as complementary hypotheses the conditions $u[1]^2 +
u[2]^2 - 2u[1]=0 \vee u[1]=0$ that can be interpreted by saying that
either  point $A$ lies on the given circle or (when $u[1]=0$) triangle
$\Delta (A, C, D)$ is rectangle at $C$.  In the next step of the
discovery procedure we consider as new hypotheses ideal the set $H +
H'$, which is of dimension 1 and where both $u[2]$ or $u[1]$ can be
taken as independent variables ruling the new construction.

\begin{verbatim}
Dim(R/(H+H'));
1
-------------------------------
Elim([x[1],x[2],x[3],x[4],u[1]],H+H');
Ideal(0)
-------------------------------
Elim([x[1],x[2],x[3],x[4],u[2]],H+H');
Ideal(0)
\end{verbatim}

\noindent Choosing, for example, $u[2]$ as relevant variable,  we
check --applying the usual automatic proving scheme-- that the new
statement $H \wedge H' \rightarrow T$ is correct under the
non-degeneracy condition $u[2] \neq 0$:
\begin{verbatim}
H'':=Elim([x[1],x[2],x[3],x[4],u[1]], Saturation(H+H',T));
H'';
Ideal(u[2]^3)
--------------------------------
\end{verbatim}
Thus we have arrived to the following statement: Given a circle of
radius 1 and centered at $(1,0)$, and a point $A$ not in the
$X$-axis and lying either on the $Y$ axis or in the circle, it
holds that the segments $AE, AF$ (where $E$ is a tangency point
from $A$ to the circle and $F$ is the intersection of the lines
passing by $(2, 0), E$  and $A, (0,0)$) are of equal length.

Let us now review Example \ref{example2} using {\it MCCGS}. As
above, the hypotheses are the union of $H:=H_1 \cup S$, where
$H_1$ expresses the equality type constraints:
\[ \begin{array}{lcl}
H_{1}&=&[(x_1-1) (u_1-x_1)+x_2 (u_2-x_2), (x_1-1)^2+x_2^2-1,\\&&
u_1 x_4-u_2 x_3, x_3 x_2-x_4 x_1-2 x_2+2 x_4]\end{array}\]
 to which we have to add the saturation ideal expressing the
inequality constraints:
\[
\begin{array}{lcl}
S &=& [u_1 x_4-u_2 x_3, x_1 u_1-u_1-x_1+x_2 u_2, x_4 x_2-2 x_2 u_2-x_3 u_1+2 u_1, \\
&& x_4 x_1-2 x_1 u_2+u_2 x_3, x_3 x_2-2 x_1 u_2+u_2 x_3-2 x_2+2 x_4, \\
&& x_1 x_3+x_3 u_1+2 x_2 u_2-2 x_1-2 u_1, x_1^2-2 x_1+x_2^2, \\
&& x_3 u_1^2+2 x_2 u_2 u_1-2 u_2^2 x_1+u_2^2 x_3-2 u_1^2-2 x_2 u_2+2 u_2 x_4, \\
&& x_3^2 u_1+x_4 u_2 x_3+2 x_4^2-4 x_3 u_1-4 u_2 x_4+4 u_1, \\
&& u_1 x_2^2-x_1 x_2 u_2-x_2^2+x_2 u_2+x_1-u_1, \\
&& u_2 x_3^3+u_2 x_4^2 x_3+2 x_4^3-4 u_2 x_3^2-4 u_2 x_4^2+4 u_2
x_3].
\end{array}
\]
The thesis is
\[ T= (u_1-x_1)^2+(u_2-x_2)^2-(u_1-x_3)^2-(u_2-x_4)^2. \]
Calling now $ \hbox{mccgs}(H_{1} \cup S \cup
T,\lex(x_1,x_2,x_3,x_4),\lex(u_1,u_2)) $ one obtains the following
segments:

\begin{center}
\begin{tabular}{||c|l|l||}
 \hline \hline
 Segment & $\lpp$ &  Description of $S_i$ \\
 \hline
 1 & $[1]$ & $K^2 \setminus (\V(u_1^2+u_2^2-2 u_1) \cup \V(u_1)) $ \\
 \hline
 2 & $[x_4^2,x_3,x_2,x_1]$ & $\V(u_1^2+u_2^2-2 u_1) \setminus (\V(u_1-2,u_2) \cup \V(u_1,u_2))$ \\
 \hline
 3 & $[x_4^2,x_3,x_2,x_1]$ & $\V(u_1) \setminus (\V(u_1,u_2^2+1) \cup \V(u_1,u_2)) $ \\
 \hline
 4 & $[x_4,x_3,x_2,x_1]$ & $\V(u_1,u_2^2+1)$ \\
 \hline
 5 & $[x_4^2,x_3,x_2^2,x_1]$ & $\V(u_1-2,u_2) $ \\
 \hline
 6 & $[x_4^2,x_3^2,x_2,x_1]$ & $ \V(u_2,u_1)$ \\
 \hline\hline
\end{tabular}
\end{center}

Segment $S_1$ states that point $A(u_1,u_2)$ must lie either in
the $Y$-axis or on the circle, as a necessary condition in the
parameter space ${\bf u}=(u_1, u_2)$ for the existence of
solutions,  in the hypotheses plus thesis variety, lying over
${\bf u}$.  This essentially agrees with the result obtained in
\cite{DR}.

A detailed analysis of the remaining segments show a variety of
formulas for determining the (sometimes not unique) values of
points $E(x_1,x_2)$ and $F(x_3,x_4)$ --verifying the theorem--
over the corresponding parameter values.

For completeness we give the different bases associated, in the
different segments,  to the above ideal of thesis plus hypotheses
\[
\begin{array}{lcl}
 B_1  &=&  [1]    \\
 B_2  &=&   [u_2^2+x_4^2-2 u_2 x_4, -u_1 x_4+u_2 x_3, u_2^3-2 u_2 u_1+x_2 u_2^2+(-2 u_2^2+2 u_1) x_4, \\
      &&  u_2 u_1+x_1 u_2-2 u_1 x_4]   \\
 B_3  &=&   [-2 u_2 x_4+x_4^2, x_3, (u_2^2+1) x_2-x_4, (u_2^2+1) x_1-u_2 x_4]   \\
 B_4  &=&  [x_4, x_3, x_2, x_1]    \\
 B_5  &=&  [x_4, -4+2 x_3, x_2^2, -2+x_1]    \\
 B_6  &=&   [x_4^2, x_3^2, -x_3 x_4+2 x_2-2 x_4, 2 x_1]   \\
\end{array}
\]
\end{example}

\begin{example}\label{Manel}

Next we consider the problem\footnote{We thankfully acknowledge here that this
problem was suggested by a colleague,  Manel Udina} described in Figure
\ref{MUdina}. Take a circle
$\cal C$ with center at $O(0,0)$ and radius 1 and let us denote
points $A=(-1, 0)$ and $B=(0,1)$. Let $D$ be  an arbitrary point
with coordinates $D=(1+a, b)$ and let $C=(1+a,0)$ be another point
in the $X$-axis, lying under point $D$. Then trace the line $BC$.
Assume this line intersects the circle $\cal C$ at point $P(x,y)$.

\begin{figure}
\begin{center}
\includegraphics{autodemo.1}

\caption{\label{MUdina} Example \ref{Manel}.}
\end{center}
\end{figure}

Consider now the, in general false, statement {\it ``the points
$A,P,D$ are aligned"}. We want to discover the conditions on the
parameters $a,b$ for the statement to be true.
The set of hypotheses plus thesis equations are very simple:
\[\begin{array}{lcl}
HT&=& [x^2+y^2-1, -x+1-y+a-a y, -2 y+b+x b-a y]
\end{array}\]
Take $x,y$ as variables and $a,b$ as parameters and call
$\hbox{\rm mccgs}(HT,\lex(x,y),\lex(a,b))$. The graphical output
of the algorithm can be seen in Figure \ref{UdinaT}, and the
algebraic description appears in the following table.

\begin{center}
\begin{tabular}{||l|l|l||}
 \hline \hline
 $\lpp$ & Basis $B_i$ & Description of $S_i$ \\
\hline $[1]$ & $[1]$   & $(K^2 \setminus (\V(a-b) \setminus \V(a-b,(b+1)^2+1))) $ \\
&  & $\cup \ (K^2 \setminus \V(2+a)) $ \\
&  & $\cup \ (K^2 \setminus \V(a-b+2))$ \\
\hline
$[y,x]$ & $[x^2+y^2-1,$   & $(\V(a-b) \setminus \V(a-b,(b+1)^2+1)) $ \\
& $x+(a+1)(y-1),$ & $\cup \ (\V(2+a) \setminus \V(b,2+a))$ \\
& $b(x+1)-(a+2)y]$ & $\cup \ (\V(a-b+2) \setminus \V(b,2+a))$  \\
\hline $[y^2,x]$ & $[y(y-1), 1+x-y]$ & $\V(b,2+a)$\\

 \hline\hline
\end{tabular}

\begin{figure}
\begin{center}
\includegraphics{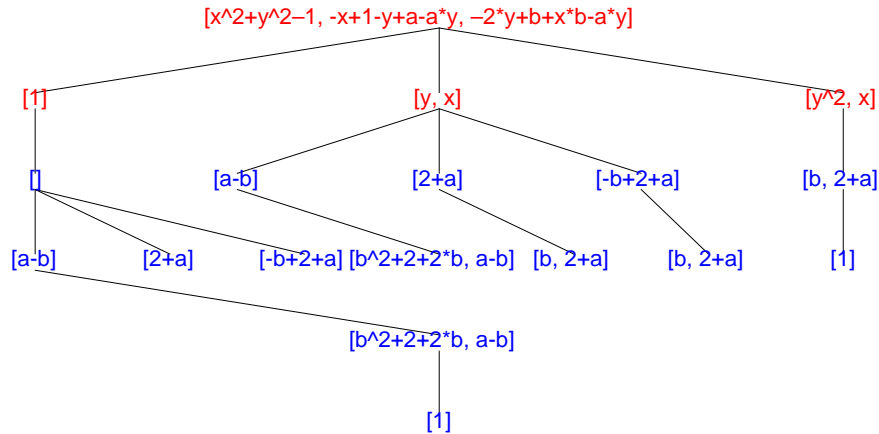}

\caption{\label{UdinaT} Canonical tree for Example \ref{Manel}}
\end{center}
\end{figure}

\end{center}

As we see, the generic case has basis $[1]$ showing that  the
statement is false in general. The interesting case corresponds,
as it is usually expected, to the case with $\lpp=[x,y]$,
providing a unique solution for $P$. The description of the
parameter set associated to this basis gives the union of three
different locally closed sets,  namely $\V(a-b) \setminus
\V(a-b,(b+1)^2+1)$, $\V(2+a) \setminus \V(b,2+a) $ and $\V(a-b+2)
\setminus \V(b,2+a)$, expressing complementary hypotheses for the
statement to hold.

The first set is (perhaps) the expected one,  corresponding to the
case $a=b$ (except for the degenerate complex point $(b,b)$ with
$(b+1)^2+1=0$, without interest from the real point of view). Thus
we can say that the statement holds if point $C$ is equidistant
from point $D$ and point $E$.

The second set yields $a=-2$ and corresponds to the situation where
 point $D$ is on the tangent to the circle trough the point
$(-1,0)$ (except for the degenerate case $b=0$). In this case
$P=A$ and, obviously,  $A, P, D$ are aligned (even in the
degenerate case, as stated in the third segment, corresponding to
the $\lpp=[y^2,x]$).

Finally, the third set gives the condition $b=a+2$ and it is also
interesting, since it corresponds to the case where the
intersecting point of the line $BC$ with the circle is taken to be
$B$ instead of $P$, and then point $D'$ should be in the vertical
of $C$ and at distance $D'C$ equal to distance $EC$ plus two.

\end{example}

\begin{example}\label{isoexample}[Isosceles orthic triangle]

In \cite{DR}  the conditions for the orthic triangle of a given
triangle (that is, the triangle built up by the feet of the
altitudes of the given triangle over each side) to the equilateral
have been discovered. Next example aims to discover  conditions
for a given triangle in order to have an isosceles orthic
triangle.

Consider the triangle of Figure \ref{isotrianglefig} with vertices
$A(-1,0)$, $B(1,0)$ and $C(a,b)$, corresponding to a generic
triangle having one side of length $2$. Denote by $P_1(a,0)$,
$P_2(x_2,y_2)$, $P_3(x_3,y_3)$ the feet of the altitudes of the
given triangle, ie. the vertices of the orthic triangle. The
equations defining these vertices are:
\[
\left.\begin{array}{lcl}
 H&=&(a-1)\,y_2-b\,(x_2-1)=0,\\ &&(a-1)\,(x_2+1)+b\,y_2=0,\\
  &&(a+1)\,y_3-b\,(x_3+1)=0,\\ &&(a+1)\,(x_3-1)+b\,y_3=0,\\
\end{array}\right\}
\]
Now let us add the condition
$\overline{P_1P_3}=\overline{P_1P_2}$.
\[
 T =(x_3-a)^2+y_3^2-(x_2-a)^2-y_2^2=0.
\]
Take $x_2,x_3,y_2,y_3$ as variables and $a,b$ as free parameters
and call
\[ \hbox{\rm mccgs}(H \cup T,\lex(x_2,x_3,y_2,y_3),\lex(a,b)).\]
The output has now four segments. The generic case, with
$\lpp=[1]$, meaning that the orthic triangle is, in general, not
isosceles; one interesting case with $\lpp=[y_3,y_2,x_3,x_2]$; and
two more  cases we can call degenerate,  with $lpp$'s
$[y_2,x_3^2,x_2]$ and $[y_2,x_3,x_2^2]$,  respectively.  For the
interesting case we show the graphic output in Figure
\ref{isotree}. Its basis is
\begin{figure}
\begin{center}
\includegraphics{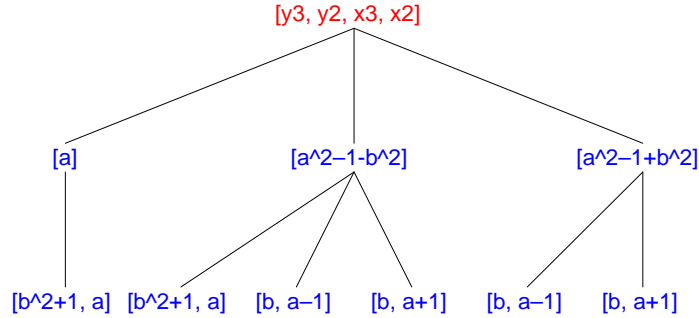}

\caption{\label{isotree} Canonical tree branch for
$\lpp=[y_3,y_2,x_3,x_2]$ in Example \ref{isoexample}.}
\end{center}
\end{figure}
 \[
 \begin{array}{lcl}
 B_2 &=& [(a^2+b^2+2a+1)y_3 -2ab-2 b,
(a^2+b^2-2a+1)y_2 +2ab-2b, \\
&& (a^2+b^2+2a+1)x_3-a^2+b^2-2a-1,
(a^2+b^2-2a+1)x_2+a^2-b^2-2a+1].
\end{array}
\]
Next table shows the description of the $lpp$ and the $S_i$'s for
the the four cases:

\begin{center}
\begin{tabular}{||l|l||}
 \hline \hline
 $\lpp$ &  Description of $S_i$ \\
\hline $[1]$ &  $K^2 \setminus ((\V(a) \setminus \V(b^2+1,a))$ \\
& $\cup\ (\V(a^2-b^2-1)\setminus \V(b^2+1,a))$ \\
& $\cup\ \V(a^2+b^2-1))$ \\
\hline $[y_3,y_2,x_3,x_2]$ &  $\V(a) \setminus \V(b^2+1,a)$ \\
& $\cup\ \V(a^2+b^2-1) \setminus (\V(b,a-1) \cup \V(b,a+1))$ \\
& $\cup\ (\V(a^2-b^2-1)\setminus (\V(b^2+1,a) \cup \V(b,a-1) \cup \V(b,a+1))$ \\
\hline $[y_2,x_3^2,x_2]$ &  $\V(b,a+1)$ \\
\hline $[y_2,x_3,x_2^2]$ &  $\V(b,a-1)$ \\
\hline\hline
\end{tabular}
\end{center}
The description of the parameter set (over the reals) for which the theorem is potentially
true and no degenerate can be phrased as follows:
\[
\begin{array}{ll}
 1)\ \ a=0 &  \\
 2)\ \ a^2+b^2=1 & \hbox{ except the points } (1,0) \hbox{ and } (-1,0)\\
 3)\ \ a^2-b^2=1 & \hbox{ except the points } (1,0) \hbox{ and } (-1,0) \\
\end{array}
\]
This set is represented in Figure \ref{fig4}.
\begin{figure}
\begin{center}
\includegraphics{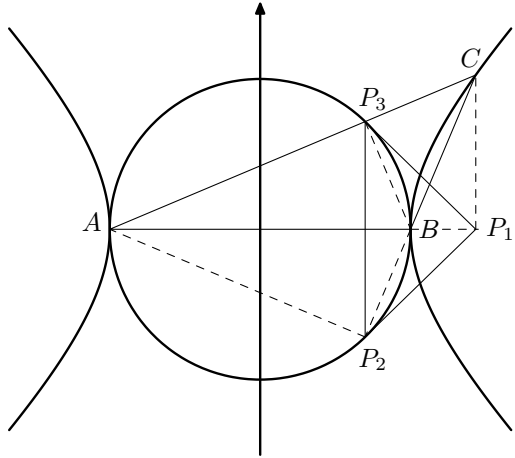}

\end{center}
\caption{\label{fig4} Solutions of Example \ref{isoexample}}
\end{figure}
and corresponds to
\begin{itemize}
\item[1)] The given triangle is itself isosceles  ($a=0$);
\item[2)] The given triangle is rectangular at vertex $C$ (with vertices $A(-1,0)$, $B(1,0)$ and
the vertex $C(a,b)$ inscribed in the circle $a^2+b^2=1$,
\item[3)] The given triangle has vertices $A(-1,0)$, $B(1,0)$ and
 vertex $C(a,b)$ lies on the hyperbola $a^2-b^2=1$.
\end{itemize}
Solution 1) is, perhaps,  not surprising. Solution 2) corresponds to rectangular
triangles for which the orthic triangle reduces to a line,
 that can be considered a degenerate isosceles
triangle. But solution 3) is a nice novelty: it exists a one parameter
family of non-isosceles triangles having isosceles orthic
triangles.

The remaining two cases in the {\it MCCGS} output with
$\lpp=[y_2,x_3^2,x_2]$ and $\lpp=[y_2,x_3,x_2^2]$  represent
degenerate triangles without geometric interest (namely $C=A$ and
$C=B$).

Thus, after performing an automatic proving procedure for the new hypotheses,
we can formulate the following theorem:
\begin{theorem} Given a triangle with vertices $A(-1,0)$,
$B(1,0)$ and $C(a,b)$, its orthic triangle will be isosceles if
and only if vertex $C$ lies either on the line $a=0$ (and then the
given triangle is itself isosceles) or in the circle $a^2+b^2=1$
(and then it is rectangular) or in the hyperbola $a^2-b^2=1$.
\end{theorem}
\end{example}

\begin{example}\label{ExSkaters}[Skaters]

Our final example is taken from the pastimes  section of the
French journal {\it Le Monde}, published on the printed edition of
Jan. 8, 2007. This example is there attributed to E. Busser and G.
Cohen. We think it is nice from {\it Le Monde} to include the
proof of a theorem as a pastime. Actually, the statement to be
proved was presented as arising from  a more down-to-earth
situation:  two ice-skaters are moving forming two intersecting
circles, at same speed and with the same sense of rotation. They
both depart from one of the points of intersection of the two
circles. Then the journal asked to show that the two skaters were
always aligned with the other point of intersection (where some
young lady, both skaters were interested at, was placed...).

Let us translate
this problem into a theorem discovering question, as follows.

We will consider two circles with centers at $P(a,1)$ and
$Q(-b,1)$ and radius $r_1^2=a^2+1$ and $r_2^2=b^2+1$, as shown in
Figure \ref{FigSkaters}, intersecting at points $O(0,0)$ and
$M(0,2)$. Consider generic points --the skaters-- $A(x_1,y_1)$ and
$B(x_2,y_2)$ on the respective circles.  Point $A$ will be
parametrized by the oriented angle $v=\widehat{OPA}$ and,
correspondingly, point $B$ will describe the oriented angle
$w=\widehat{OQB}$. Therefore we can say that angle zero
corresponds to the departing location of both skaters, namely,
point $O$.

We claim that, for whatever position of points $A, B$, {\it the
points $A,M,B$ are aligned}, which is obviously false in general.
But we want to determine if there is a relation between the two
oriented angles making this statement to hold true. Denote
$c_v,s_v,c_w,s_w$ the cosine and sine of the angles $v$ and $w$.
It is easy to establish the basic hypotheses, using scalar
products:

\[\begin{array}{lcl} H_1&=&[(x_1-a)^2+(y_1-1)^2-a^2-1,(x_2+b)^2+(y_2-1)^2-b^2-1, \\
&& a (x_1-a)+(y_1-1)+(a^2+1) c_{v}, -b (x_2+b)+(y_2-1)+(1+b^2)
c_{w}
\end{array}\]
\begin{figure}
\begin{center}
\includegraphics{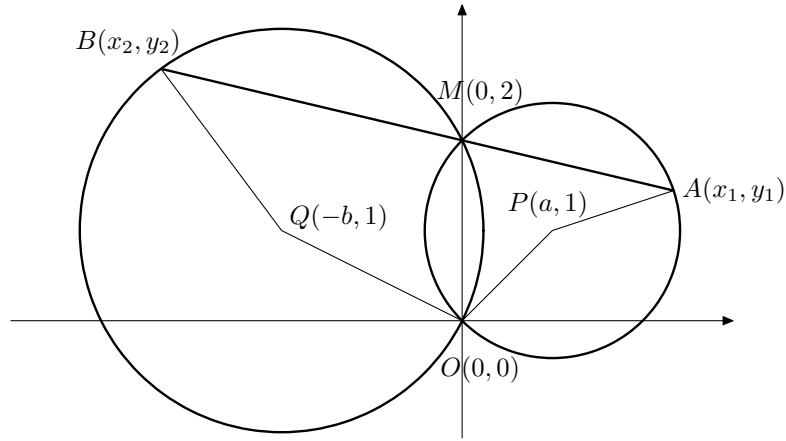}

\caption{\label{FigSkaters} Skaters problem}
\end{center}
\end{figure}
Now, as the angles are to be taken oriented (because we assume the
skaters tare moving on the corresponding circle in the same
sense),  we need to add the vectorial products involving also the
sine to determine exactly the angles and not only their cosines.
So we add the hypotheses:
\[\begin{array}{lcl}
H_2 &=& [a (y_1-1)-(x_1-a)+(a^2+1) s_v,-b (y_2-1)-(x_2+b)+(b^2+1) s_w]\\
\end{array}\]
The thesis is,  clearly:
\[T= x_1 y_2-2 x_1-x_2 y_1+2 x_2. \]

The radii of the circles are
\[r_1^2 = a^2+1 \ \ \hbox{\rm and } r_2^2 = b^2+1\]
and for $r_1 \ne 0$ and $r_2 \ne 0$ we have
\[
\begin{array}{ll}
c_{v_0} =\cos v_0 = \cos \widehat{OPM} = \displaystyle\frac{a^2-1}{a^2+1},
&  s_{v_0}=\sin v_0 = \sin \widehat{OPM} = \displaystyle\frac{-2a}{a^2+1},\\
c_{w_0} =\cos w_0= \cos \widehat{OQM}=\displaystyle\frac{b^2-1}{b^2+1},
& s_{w_0}=\sin w_0 = \sin \widehat{OQM}=\displaystyle\frac{2b}{b^2+1}.\\
\end{array}
\]
We want to take $a,b$ and the angles $v$ and $w$ --in terms of the
sines and cosines-- as parameters. So we must introduce the
constraints on the sine and cosine parameters. Moreover, we notice
there are also some obvious degenerate situations, namely $r_1=0$,
$r_2=0$ and $a+b=0$,  corresponding to null radii or coincident
circles, and we want to avoid them.

Currently,  {\it MCCGS} allows us to introduce all these
constraints in order to discuss the parametric system. The call is
now
\[
\begin{array}{l}
\hbox{\rm mccgs}(H_1 \cup H_2 \cup
T,\ \lex(x_1,y_1,x_2,y_2),\ \lex(a,b,s_v,c_v,s_w,c_w),\\
\hspace{1cm}  \hbox{null}=[c_v^2+s_v^2-1,c_w^2+s_w^2-1],\
\hbox{notnull}=\{a^2+1,b^2+1,a+b\}). \end{array}
\]
including the constraints on the parameters and eluding degenerate
situations as options for {\it MCCGS}.

The result is that {\it MCCGS} outputs only 2 cases. The first one
has basis $[1]$,  showing that, in general, there is no solution
to our query. The second one has $\lpp=[y_2,x_2,y_1,x_1]$
determining in a unique form the points $A$ and $B$ for the given
values of the parameters. The associated basis is
\[[y_2+c_w-b s_w-1,x_2-b c_w-s_w+b,y_1+c_v+a s_v-1,x_1+a c_v-s_v-a]\]
with parameter conditions that can be expressed as the union of
three irreducible varieties:
\[ \begin{array}{lcl}
V_1 &=& \V(c_w^2+s_w^2-1, c_v-c_w, s_v-s_w) \\
V_2 &=& \V(c_w^2+s_w^2-1, c_v^2+s_v^2-1, s_w+b c_w-b,b s_w-c_w-1)
\\
V_3 &=& \V(c_w^2+s_w^2-1, c_v^2+s_v^2-1, -s_v+a c_v -a, a
s_v+c_v+1)\\
\end{array}
\]

The interpretation is easy: $V_1$ corresponds to arbitrary $a, b,
w$,  plus the essential condition $v=w$, which is the interesting
case, stating that our conjecture requires (and it is easy to show
that this condition is sufficient) that both skaters keep moving
with the same angular speed.

$V_2$ corresponds to $s_w=s_{w_0},\ c_w=c_{w_0}$ and $a, b, v$ free,
thus $B=M$ and $A$ can take any position.

$V_3$ is analogous to $V_2$,  and corresponds to placing $A=M$ and $B$ anywhere.

So we can summarize the above discussion in the following
 \begin{theorem} Given two non coincident
circles of non-null radii and centers $P$ and $Q$, intersecting at
two points $O$ and $M$, let us consider  points $A$, $B$ on each
of the circles. Then the three points $A,M,B$ are aligned if and
only if the oriented angles $\widehat {OPA}$ and $\widehat{OQB}$
are equal or
  $A$ or $B$ or both coincide with
 $M$.
 \end{theorem}
 \end{example}

 \section{Performances}

Although the principal advantage of {\it MCCGS} in relation to
other CGS algorithms is the simplicity and properties of the
output: the minimal number of segments and the characterization of
the type of the solution depending on the values of the
parameters, the computer implementation\footnote{That can be freely
obtained at http://www-ma2.upc.edu/$\sim$montes}  of  the corresponding package, named  {\it dpgb}
release 7.0,  in {\it Maple 8} is relatively short time consuming.
Moreover,we think that no other actual PCAD software will be able to obtain
the accurate result obtained, for example, in example 13.
We give here a table with the CPU time and number of segments for
the examples of the paper.\vspace{0.2cm}

\begin{center}
\begin{tabular}{||c|c|c||}
\hline\hline
Example & CPU time (sec.) & Number of segments \\
\hline\hline
\ref{singpointsconic} & 1.9 & 3 \\
\hline \ref{example2} & 12.8 & 6 \\
\hline \ref{Manel} & 0.98 & 3\\
\hline \ref{isoexample} & 4.4 & 4\\
\hline \ref{ExSkaters} & 129.4 & 2\\
 \hline\hline
\end{tabular}
\end{center}\vspace{0.2cm}

The computations were done with a Pentium(R) 4 CPU at 3.40 Ghz and
1.00 GB RAM.

\section{Conclusion}
We have briefly introduced the principles of automatic discovery
and also the ideas --in the context of comprehensive Gr\"obner
basis-- for discussing polynomial systems with parameters,  via
the new  {\it MCCGS} algorithm. Then we have shown how natural is
to merge both concepts, since the parameter discussion can be
interpreted as yielding, in particular, the projection of the
system solution set  over the parameter space; and since the
conditions for discovery can be obtained by the elimination of the
dependent variables over the ideal of hypotheses and thesis. Moreover,
we have also remarked how the approach through {\it MCCGS} provides
new candidate complementary conditions of more general type and, in
some particular instances (segments of the parameter space yielding
to unique solution),  quite common in our examples, an easy test for
the sufficiency of these conditions. Finally, the use of {\it MCCGS}
for automatic proving has been presented, as part of a  formal discussion
on the limitations of the discovery method.

We have exemplified this approach through a collection of
non-trivial examples (performed by running the current Maple
implementation of {\it MCCGS}, see \cite {MaMo06},  over a
laptop, without special time -- a few seconds--  or memory
requirements), showing that in all cases, the {\it MCCGS}  output
is very suitable  to providing geometric insight, allowing the
actual discovery of interesting and new? theorems (and pastimes!).



\end{document}